\magnification=\magstep1
\advance\voffset by 1 true cm
\baselineskip=15pt
\overfullrule =0pt 
\font\bbb=msbm10

\def\C{\hbox{\bbb C}}
\def\P{\hbox{\bbb P}}

\centerline {\bf On the connectivity of complex affine hypersurfaces, II}

\vskip2truecm

\centerline {\bf by  Alexandru Dimca and Laurentiu Paunescu}

\vskip2truecm

Let $f \in \C[x_0,x_1, ... ,x_n]$ be a polynomial and $ {\bf w}=(w_0,...,w_n)$
 a system of positive integer weights.
Consider $f=f_0+ ...+f_e+f_d$  the decomposition of $f$ as a sum of weighted
homogeneous components $f_i$ where $deg(f_i)=i$ with respect to ${\bf w}$. Here $f_d$ is the top degree nonzero component and we assume that $f_j=0$ for $e<j<d$ for some integer $e>0$.

For any polynomial $g$ we denote by $\partial g$ its gradient

$$ \partial g =( g_{x_0},...,g_{x_n})$$
where $g_{x_i}$ is the partial derivative of $g$ with respect to $x_i$.

Define the subset $S(f)$ in $\C^{n+1}$ given by the equations $\partial f_d=0, f_e=0$. With this notation our result is the following.

\bigskip

\noindent {\bf THEOREM}

(i) Let $\Sigma(f)=\{x \in \C^{n+1}; \partial f =0\}$ be the set of singular points of the polynomial $f$. Then dim$\Sigma(f) \leq $ dim$S(f)$.

(ii) Any fiber of the polynomial $f$ is $(n-1-$dim$S(f))$-connected.

\bigskip

Improving previous results on the connectivity of affine hypersurfaces, the first author has proved the part (ii) of this result in the case of the usual weights ${\bf w}=(1,1,...,1)$ in [D1]. 

The proof there was based on the properties of tame and quasi-tame polynomials introduced by Broughton [B],
N\'emethi [N1], [N2] and N\'emethi-Zaharia [NZ]. More precisely, the general case was reduced to these special classes of polynomials by applying an affine  Lefschetz hyperplane section theorem due to Hamm [H].
This proof cannot be applied for arbitrary weights since in general we don't have enough weighted hyperplanes to proceed by induction (even though Lemma 9. in [D1] can be easily extended to show that the polynomial $f$ is quasitame when dim$S(f)=0$). 

On the other hand, N\'emethi and Sabbah [NS] have recently studied tame polynomials defined on affine varieties, but their definition does not seem easily adaptable to keep track of the chosen weights.

Our proof below is based on a powerful general result by Hamm, namely Prop.3
in [H] (which is the key step in his proof of the affine Lefschetz theorem on hyperplane sections and which we slightly correct) and on basic properties of isolated complete intersection singularities, for which we refer the reader to Looijenga's book [L].

\bigskip

The following example shows the strength of the new result over the old.

\bigskip

\noindent{\bf Example 1.}

Let $f:\C^6 \to \C$ be the polynomial given by 
$$f(x)= x_0^3x_1^3x_2^3+(x_0+x_1+x_2)^7+x_3^5+x_4^4+x_5^3+(x_0+x_1+x_2+x_3+x_4+x_5)^2$$ 
If we consider the usual weights ${\bf w}=(1,1,...,1)$, then dim$S(f)=4$, and hence we get that the fibers of $f$ are 0-connected, i.e. they are connected. In this case $d=9$, $e=7$

Now, if we use the system of weights ${\bf w}=(20,20,20,36,45,60)$ chosen such that the top degree form $f_d$ contains as many monomials as possible, then
dim$S(f)=1$ hence the fibers of $f$ are in fact 3-connected. In this case $d=180$, $e=140$.

\bigskip

It is interesting to notice that the topological result above has some useful
algebraic consequences expressed in terms of various complexes of differential forms associated naturally to the polynomial $f$.

To state them, let $A^*=(\Omega ^*, $d$)$ denote the De Rham complex of global regular differential forms on $\C^{n+1}$ with d the exterior differentiation acting on forms (not to be confused to the degree $d$ which occured above). 

The first complex associated to $f$ is the complex $K_f^*=(\Omega ^*, $d$f \wedge)$ which  can be identified to the Koszul complex of the partial derivatives of $f$ in the polynomial ring $\C[x_0, ... x_n]=\Omega ^0$.

The De Rham complex $A^*$ has a natural subcomplex $B_f=($d$f \wedge \Omega ^*, $d$)$ and a natural quotient complex $C_f=A^*/B_f^*$ called the complex of global
relative differential forms.

Finally, one can consider as in Dimca-Saito [DS], the  mixture of De Rham and Koszul complexes, namely the complex $D_f=(\Omega ^*,$ d$-$d$f \wedge)$.

\bigskip

\noindent {\bf Corollary.}

Let $f$ be a polynomial as above and assume that $k=$dim$S(f)<n$. Then

(i) $H^i(K_f^*)=0$ for all $0 \leq i \leq n-k$;

(ii) $H^1(B_f^*)=\C[f]$d$f$ and $H^i(B_f^*)=0$ for all $ i \not=1, 0 \leq i \leq n-k$;

(iii) $H^0(C_f^*)=\C[f]$ and  $H^i(C_f^*)=0$ for all $0 < i \leq n-k-1$;

(iv) $H^i(D_f^*)=0$ for all $0 \leq i \leq n-k$.

\bigskip

{\it Proof of the Theorem}

\bigskip

Let $N$ be any common multiple of the weights $w_i$ and consider the integers
$m_i=N/w_i$. Let $\phi: \C ^{n+1} \to \C^{2n+2}$ be the embedding

$$ \phi(x)=(x_0,...,x_n,x_0^{m_0},...,x_n^{m_n})$$

Let $X=f^{-1}(0)$. It is enough to prove that

\bigskip

(i') dim$X_{sing} \leq $ dim$S(f)$;

(ii') $X$ is $(n-1-$dim$S(f))$ connected.

\bigskip

Let $Y=\phi(X)$ and note that $\phi:X \to Y$ is an isomorphism.
Let $k=$dim$S(f)$ and let $H_1:l_1=0$,...,$H_k:l_k=0$ be generic hyperplanes in $\C^{2n+2}$.

Using repeatedly Lefschetz hyperplane section theorem in [H] we get that the
inclusion $ Y_0=Y \cap H_1 ...\cap H_k \to Y$ is an $(n-k)$-equivalence.
 
Let $X_0=\phi ^{-1}(Y_0)$ and $g_i(x)=l_i(\phi (x))$ be the corresponding
pull-back polynomials on $\C^{n+1}$. Note that the weighted homogeneous component $g_{i,N}$ is a polynomial of Pham-Brieskorn type, i.e. a sum $\sum a_{ij}x_j^{m_j}$. 

It is easy to see by general transversality arguments involving the parameter space of all the coefficients $a_{ij}$ that we can choose the generic hyperplanes $H_i$ such that one has the following properties:

\bigskip

(a) $A= \{ x \in \C^{n+1}; g_1(x)=...=g_k(x)=0 \}$ is a smooth complete intersection of dimension $n+1-k$;

\bigskip

(b) the set $\{x \in \C^{n+1}; x \in S(f), g_{1,N}(x)=...=g_{k,N}(x)=0\}$ is just the origin;

\bigskip

(c) there exist weighted homogeneous polynomials $h_1,...,h_n$ such that the
germs at the origin of the sets $\{x \in \C^{n+1}; g_{1,N}(x)=...=g_{k,N}(x)=f_d(x)=h_1(x)=...=h_j(x)=0 \}$ is an isolated complete intersection singularity
of dimension max$(n-k-j,0)$ for all $j=1,...,n$.

\bigskip

Note also that we can arrange such that
 
$(Y_0)_{sing}=Y_{sing}\cap H_1 \cap ...\cap H_k$ and  dim$(Y_0)_{sing}$=dim$Y_{sing}-k$ 

In particular, (i') will follow
from the next result.

\bigskip

\noindent{\bf Lemma 1.}

$X_0$ has at most isolated singularities.

\bigskip

{\it Proof of Lemma 1.}

If $X_0$ has not isolated singularities, then an application of the curve selection lemma, see Milnor [M], would produce a path $p: (0, \epsilon) \to (X_0)_{sing}$ given by a Laurent power series
$$p(t)=c_0t^s+c_1t^{s+1}+...$$
where $c_i \in \C^{n+1}$, $c_0 \not= 0$ and the integer $s$ is strictly negative.

Note that we can reparametrize this path by replacing $t$ with $t^N$, i. e. we
can suppose that all exponents in $p(t)$ with non zero coefficients are divisible by our weights $w_i$.

There is a standard $\C^*$-action on $\C^{n+1}$ associated with these weights,
namely
$$ u*x=(u^{w_0}x_0,...,u^{w_n}x_n)$$
We can rewrite our path using this action in the form
$$ p(t)=t^s *c_0+t^{s+1}*c_1+...$$
where the coefficients $c_i$'s and the exponent $s$ are different in general from those considered first above, but which enjoy the same properties, i.e.
$c_0 \not= 0$ and $s<0$.

Next we have
$g_i(p(t))=g_{i,N}(c_0)t^{sN}+...=0$ where the dots represent higher order terms, hence we have $g_{i,N}(c_0)=0$.
In a similar way we get $\partial f_d(c_0)=0$.

We have also

$0=df(p(t))- \sum _{i=0,n} w_if_{x_i}(p(t))p_i(t)= (d-e)f_e(p(t))+...=
(d-e)f_e(c_0)t^{se}+...$ which implies that $f_e(c_0)=0$, a contradiction
with the property (b) above.

This ends the proof of the Lemma and of part (i) in our Theorem.

\bigskip

To continue the proof, we would like to apply Proposition 3. in Hamm [H], and this statement has to be slightly corrected. Indeed, with the notation from [H],
consider the following example.

\bigskip

\noindent{\bf Example 2.}

\bigskip
Let $A$ be $\C^2$ with the trivial Whitney stratification, $f_1=x^2y+x$ and
$f_2=x$. Then all the assumptions in Proposition 3. in [H] are fulfilled, but not the conclusion, i.e. $\C^2$ cannot be obtained from the non connected space
$f_1^{-1}(0)$ by adjoining 2-cells.

\bigskip

However, looking at Hamm's proof of Proposition 3., it is clear that one needs the following additional condition.
\bigskip

$(\star)$ the set $\{z \in A; |f_1(z)| \leq a_1, ..., |f_N(z)| \leq a_N \}$ is compact for any positive numbers $a_j$, $j=1,...,N$.

\bigskip

This condition allows one to use the Morse theory for the function $|f_1|^2$ to increase $a_1$ for some fixed (and very large) $a_2,...,a_N$ as in Hamm's proof. Note that this condition $(\star)$ is not fulfilled in our Example 2.

\bigskip

To prove the claim (ii') we use this completed version of Proposition 3. in [H] for the following data:

The set $A$ in [H] is here the smooth complete intersection $A=\{x \in \C^{n+1};g_1(x)=...=g_k(x)=0 \}$, see property (a) above, and hence it satisfies the necessary local connectivity conditions. We choose the trivial Whitney stratification, in which there is just one stratum, namely $A$.

The functions $f_i$ which appear in [H] are chosen as follows $f_1=f$ and
$f_j=h_{j-1}$ for $j=2,...,n+1$.

The conditions these functions have to satisfy are the following.

\bigskip

(c1) the connected components of the critical set of the function $f_1:A \to \C$ are compact.

This condition is fulfilled, since the same computations as in our Lemma above shows
that these components are actually points.

\bigskip

(cj) Let $F$ be the critical set of the function $(f_1,...,f_j):A \to \C^j$
for some $j>1$. Then it is enough to show that $(f_1,...,f_{j-1}):F \to \C ^{j-1}$ is a proper mapping.

In our case, a critical point for $(f_1,...,f_j):A \to \C^j$ is just a singular point of the variety $Z$ given by the equations

$$ g_1(x)=...=g_k(x)= 0, f(x)= a_1, h_1(x)=a_2, ..., h_{j-1}(x)=a_j$$

This variety corresponds to a fiber in a deformation of the isolated complete intersection singularity

$$g_{1,N}(x)=...=g_{k,N}(x)= f_d(x)=h_1(x)=...=h_{j-1}(x)=0$$
For details on how an affine variety can be ``localized'' to become a fiber in a singularity deformation see [D2], p. 157 and p. 161 (the case discussed there is for the usual weights ${\bf w}= (1,...,1)$ but the same idea works for general weights).

Hence any such variety $Z$ is either smooth or has just isolated singularities.
Projecting on the space $\C^{j-1}$ means that now we consider a family of fibers
as above corresponding to a line in the base of the deformation of our singularity. If $h_j$ is chosen general enough, this line is not contained in the discriminant $\Delta$ of this deformation, i.e. there is just a finite number of singular fibers when we vary $a_j$. This follows from the basic fact that a nonconstant regular function on variety with isolated singularities has just finitely many critical values.

It follows that the map $(f_1,...,f_{j-1}):F \to \C ^{j-1}$ is the composition of two finite maps, the first from $F$ to the discriminant $\Delta$ of the deformation and the second the projection $\Delta \to \C^{j-1}$.
 This explains why the condition (cj) holds as well
in our setting for $j>1$. 

Note also that when $j > n+1-k$ , then $F=A$ and the map $(f_1,...,f_{j-1}):A \to \C^{j-1}$ is a finite map. This shows in particular that our additional condition $(\star)$ holds in our case.

\bigskip

By Proposition 3 in [H] we have that the inclusion $X_0 \to A$ in an $(n-k)$-equivalence. As above, we can identify $A$ with a smooth fiber (Milnor fiber)
in the deformation of the corresponding isolated complete intersection,
hence $A$ is $(n-k)$-connected being a bouquet of $(n-k+1)$-spheres.

Combining this with the $(n-k)$-equivalence $X_0 \to X$ obtained at the beginning, we have that $X$ is $(n-k-1)$-connected.

This ends the proof of our Theorem.

\bigskip

{\it Proof of the Corollary}

The first claim (i) depends only on part (i) of our Theorem. Indeed, the cohomology groups of the Koszul complex are finitely generated $\C[x_0,...,x_n]$-modules and to prove that one of them is trivial it is enough to show that all its localizations at maximal ideals are trivial. 

To check this local property we can use GAGA and replace algebraic localization by analytic localization. At this level the result follows from a general result in Looijenga's book [L], namely Corollary (8.16), p. 157 (take $X$ a smooth germ and $k=1$ in that statement).

To prove (ii) and (iii), we consider the exact sequence of complexes
$$ 0 \to B_f^* \to A^* \to C_f^* \to 0$$
This shows that it is enough to prove (ii). We have in fact the following more precise result. By convention, the dimension of the empty set is taken to be $ -1.$

\bigskip

\noindent{\bf Lemma 2.}

Let $f:\C^{n+1} \to \C$ be a polynomial function whose critical locus $\Sigma(f)$ has dimension $k<n$. Then the following statements are equivalent:

(i) $H^1(B_f)=\C[f]$d$f$ and $H^i(B_f)=0$ for $i \not= 1, i \leq n-k$;

(ii) the reduced cohomology groups ${\tilde H}^i(F_t,\C)$ are trivial for $0 \leq i \leq n-k-1$.

\bigskip

{\it Proof of Lemma 2.}

This proof follows closely Section 9. in Sabbah [S] and for this reason we just mention the main new points.

The algebraic Gauss-Manin system of $f$ is represented (up to a shift) by the complex of $A_1=\C[t]<\partial >$-modules $GM_f^*= (\Omega ^* [\partial], d_f)$ where the $\C$ -linear differential $d_f$ is defined by $d_f(\omega \partial ^m)=$d$ \omega
\partial ^m -$d$f \wedge \omega \partial ^{m+1}$.

As in Prop. (9.2) in [S], it follows that the condition (ii) is equivalent to
the condition 
\bigskip

(ii') $H^1(GM_f^*)=\C[f]$d$f$ and $H^i(GM_f^*)=0$ for $i \not= 1, i \leq n-k$.

The complex $GM_f^*$ comes equipped with a decreasing filtration given by

$$F^sGM_f^m= \Omega ^m [ \partial ]_{\leq m-s}$$
where the filtration on the RHS is by the degree with respect to $\partial$
(Sabbah prefers to work with a similar but increasing filtration in [S]).

The general theory of spectral sequences, see if necessary [Mc], associates to this decreasing, exhaustive and bounded below filtration a spectral sequence with $E_1^{s,t}=H^{s+t}(Gr^s_FGM_f^*)$ converging to $H^{s+t}(GM_f^*)$.

For $t>0$, we have $E_1^{s,t}=H^{s+t}(K_f^*)$ and hence in particular in our case we have $E_1^{s,t}=0$ for all $t>0,s+t \leq n-k$ by our Corollary (i) above.

Moreover, the terms $E_1^{s,0}$ with the corresponding differential $d_1:E_1^{s,0} \to E_1^{s+1,0}$ coming from the spectral sequence can be identified for $s < n-k$ (since we need again Corollary (i)) to the corresponding initial part in the complex $B_f^*$.

Since this part of the spectral sequence clearly degenerates at the $E_2$-term, i.e. $E_2^{s,0}=E_{\infty}^{s,0}$ for $s \leq n-k$, we obtain the equivalence $(i) \Leftrightarrow (ii)'$ which
completes the proof of Lemma 2. 

To end the proof of our Corollary, we have just to use the main results in [DS], saying that the cohomology of the complex $D_f$ twisted by -1 is just the reduced cohomology of the general fiber of the polynomial $f$.

\bigskip
 
\noindent {\bf Remarks.}

(I) The condition $k=dim \Sigma(f)<n$ is equivalent to saying that the fibers of $f$ are all reduced. This condition is not needed for part (i) of our Corollary (which holds even when $k=n$, but it is necessary for parts (ii)-(iv) since it
implies that the general fiber of $f$ is connected.

Indeed, it is well known that any polynomial $f$ can be written as a composition $f(x)=h(g(x))$ where $g: \C^{n+1} \to \C$ has a connected general fiber and $h:\C \to \C$, both $g$ and $h$ being polynomials. Since it is difficult to locate a reference for this fact, here is a short proof of it suggested  to the first author by Zaidenberg, to whom we are grateful.

Assume that $f$ is not a constant polynomial and let ${\tilde f}: X \to \P ^1$ be a smooth compactification of $f$. Then the Stein Factorization Theorem, see [Har], p.280, gives a smooth curve $C$ and morphisms $ \tilde g:X \to C$ and $\tilde h: C \to \P ^1$ such that $ \tilde h \circ \tilde g = \tilde f$ and such that all the fibers of $\tilde g$ are connected.

A generic line $L$ in $\C^{n+1}$ has the following properties:

(i) $f$ is not constant when restricted to $L$;

(ii) the closure $\tilde L$ of $L$ in $X$ is a smooth rational curve which meets ${\tilde f}^{-1}(\infty)$ at exactly one point and this intersection is transverse.

{\noindent Then} $\tilde L$ is a rational curve and $\tilde g | \tilde L$ is a non constant map. This implies that $C= \P^1$.
Moreover, we have ${\tilde h}^{-1}(\infty)=\infty$, otherwise the condition (ii) above is contradicted.
This implies that $\tilde h$ give by restriction to $\C=\P^1 \setminus \{\infty \}$ a map $h: \C \to \C$ such that $f=h \circ g$ where $g= \tilde g | \C^{n+1}$ has its general fiber connected.

\bigskip

(II) For a given polynomial $f: \C^{n+1} \to \C$ and a given system of weights
$\bf w$ we define the connectivity order of $f$ with respect to $\bf w$ to be the integer
$$ c_{\bf w}(f)=n-1 -k$$
where $k=$dim$S(f).$

Then it is easy to see that for the sum of two polynomials $f+g$ where $g:\C^{m+1} \to \C$ and for any weights $\bf w'$ associated to $g$ such that the corresponding degrees $d$ for $f$ and $g$ coincide (by multiplying any weights $\bf w$ for $f$ and $\bf w'$ for $g$ by suitable positive integers this can be always achieved) we have
$$ c_{\bf w,w'}(f+g) \geq c_{\bf w}(f)+c_{\bf w'}(g) +1$$
in spite of the fact that in general the corresponding degrees $e$ will be different.

Note that the general fiber of $f+g$ is the join of the general fibers of $f$ and $g$, see [N3], but nothing is known about the special fibers.

\bigskip

(III) Siersma and Tib\u ar have shown in [ST] that in the case of the usual weights ${\bf w}=(1,1,...,1)$ and if $e=d-1$ and dim$\Sigma(f)<$dim$S(f)=k$, then the general fiber of $f$ is $(n-k)$-connected, i.e. a better by 1 estimate than
that given by our Theorem. They also give an example showing that this better estimate
fails for the special fibers.

Improving the results in [ST], Tib\u ar has recently obtained very general connectivity results for the fibers of polynomial mappings, see [T], Thm . (5.5). However the effectiveness of his results depends on the explicit construction of
fiber-compactifying extensions of the polynomial function $f$ having a small
critical set at infinity. It seems to us that  such constructions are difficult to handle for general weights ${\bf w}$ or when $e<d-1$.

\vskip2truecm

\noindent {\bf Acknowledgements.}

The first author would like to thank Dan Burghelea and Andr\'as N\'emethi for
their kind invitation to visit  Ohio State University at Columbus.

The second author thanks the staff of the Tokyo Institute of Technology for
their hospitality.

This research was supported in part by the ARC Grant A69530267.

\vskip3truecm

\noindent {\bf REFERENCES}
\bigskip

\item{[B]} S. A. Broughton: Milnor numbers and the topology of polynomial hypersurfaces, Invent. Math. 92(1988),217-241.

\item{[D1]} A. Dimca: On the connectivity of complex affine hypersurfaces,
Topology, 29(1990), 511-514.

\item{[D2]} A. Dimca: Singularities and Topology of Hypersurfaces,
Universitext,Springer, 1992.

\item{[DS]} A. Dimca, M. Saito: On the cohomology of the general fiber of a polynomial map, Compositio Math. 85(1993), 299-309.

\item{[H]} H. Hamm: Lefschetz theorems for singular varieties, Proc. Symp. Pure Math. 40, Part I, Amer. Math. Soc. (1983), 547-557.

\item{[Har]} R. Hartshorne: Algebraic Geometry, Graduate Texts in Maths. 52, Springer-Verlag, Berlin, 1977.

\item{[L]} E. J. N. Looijenga: Isolated Singular Points on Complete Intersections, London Math. Soc. Lecture Note Series 77, Cambridge University Press, Cambridge, 1984.

\item{[Mc]} J. McCleary: User's Guide to Spectral Sequences, Publish or Perish, 1985.

\item{[M]} J. Milnor: Singular Points of Complex Hypersurfaces, Ann.Math. Stud.
61, Princeton University Press, Princeton (1968).

\item{[N1]} A. N\'emethi: Th\'eorie de Lefschetz pour les vari\'et\'es
alg\'ebriques affines, C. R. Acad. Sci. Paris 303(1986), 567-570.

\item{[N2]} A. N\'emethi: Lefschetz theory for complex affine
varieties, Rev. Roum. Math. Pures et Appl. 33(1988),233-260.

\item{[N3]} A. N\'emethi: Generalized local and global Sebastiani-Thom type theorems, Compositio Math. 80(1991),1-14.

\item{[NS]} A. N\'emethi, C. Sabbah: Semicontinuity of the spectrum at infinity, preprint.

\item{[NZ]} A. N\'emethi, A. Zaharia: On the bifurcation set of a
polynomial function and Newton boundary, Publ.RIMS Kyoto
Univ. 26(1990), 681-689.

\item {[S]} C. Sabbah: Hypergeometric periods for a tame polynomial, preprint
1996.

\item{[ST]} D. Siersma, M. Tib\u ar: Singularities at infinity and their vanishing cycles, Duke Math. J. 80(1995), 771-783.

\item {[T]} M. Tib\u ar: Topology at infinity of polynomial mappings and Thom regularity condition, Compositio Math. 111(1998), 89-109.

\bigskip

Laboratoire de Math\'ematiques Pures de Bordeaux

Universit\'e Bordeaux I

33405 Talence Cedex, FRANCE

\bigskip

email: dimca@math.u-bordeaux.fr

\bigskip

School of Mathematics and Statistics,

Sydney University,

Sydney 2006, AUSTRALIA

\bigskip

email:laurent@maths.usyd.edu.au

\bye